\documentclass[reqno]{amsart}
\usepackage{amssymb, amsmath, amsthm, amsthm}
\usepackage[bookmarks=false]{hyperref}

\newtheorem*{rem}{Remark} 

\newtheorem*{rems}{Remarks} 

\newtheorem*{sol}{Solution}

\newtheorem*{nt}{Notes} 

\makeatletter
\renewcommand{\@biblabel}[1]{#1.}
\makeatother

\newcommand{\qrfac}[2]{{\left({#1}; q\right)_{#2}}}

\begin{document}

\title{How to Discover the Rogers--Ramunujan Identities}
\author{Gaurav Bhatnagar}
\address{ Educomp Solutions Ltd., Delhi, India}
\email{bhatnagarg@gmail.com}
\date{January 18, 2012}

\begin{abstract}
We examine a method to conjecture two very famous identities that were conjectured by Ramanujan, and later found to be known to Rogers. 
\end{abstract}

\maketitle


\section{Introduction}
The so-called Rogers--Ramunujan identities were sent by Ramanujan to Hardy nearly 100 years ago.
In the next few years, the identities were circulated amongst mathematicians, but nobody, including Ramanujan, was able to prove them. Then one day, while riffling through old back copies of a journal, Ramanujan himself discovered them in an obscure paper written in 1894 by the English mathematician Rogers.  This spurred both Rogers and Ramanujan to provide simpler proofs of the identities, that were published in 1919. 

The identities are:
\begin{align}\label{rr1-expanded}
1+&\frac{q}{(1-q)}+\frac{q^4}{(1-q)(1-q^2)}
+\frac{q^9}{(1-q)(1-q^2)(1-q^3)}+\cdots \cr
&=\frac{1}{(1-q)(1-q^6)(1-q^{11})(1-q^{16})\cdots}
\times \frac{1}{(1-q^4)(1-q^9)(1-q^{14})\cdots};\cr
\end{align}
and,
\begin{align}\label{rr2-expanded}
1+&\frac{q^2}{(1-q)}+\frac{q^6}{(1-q)(1-q^2)}
+\frac{q^{12}}{(1-q)(1-q^2)(1-q^3)}+\cdots \cr
&=\frac{1}{(1-q^2)(1-q^7)(1-q^{12})(1-q^{17})\cdots }
\times
\frac{1}{(1-q^3)(1-q^8)(1-q^{13})\cdots }.\cr
\end{align}
About these Hardy  \cite[p.~xxxiv]{hardy} remarked:
\begin{quote}
It would be difficult to find more beautiful formulae than the \lq\lq Rogers--Ramunujan" identities \dots
\end{quote}

The purpose of this article is to introduce you to the Rogers--Ramunujan identities, by discussing an approach to discover them. When you see that they appear from a very simple generalization of the simplest possible infinite continued fraction, that in turn is related to the celebrated Fibonacci sequence, perhaps you may begin to agree with Hardy's opinion of these formulas.

\section{On infinite series and products}

Examine the first of the Rogers--Ramunujan identities \eqref{rr1-expanded}. On the left hand side, there is an infinite sum of terms such as
$$\frac{q^{k^2}}{(1-q)(1-q^2)\cdots (1-q^k)}.$$
On the other side, there is an infinite product, that has terms such as 
$$\frac{1}{1-q^m}$$
where $m$ is either of the form $5k+1$ or $5k+4$, for $k=0, 1, 2, \dots$. 
Our first goal is to understand how to assign meaning to such infinite sums and products. 
%
%

To that end, let us examine the simplest such identity, namely  the {\bf geometric series}:
\begin{equation}\label{geometric-series}
1+q+q^2+q^3+\cdots = \frac{1}{1-q}, \text{ where $|q|<1$}
\end{equation}
Here the left hand side is an infinite sum of terms such as $q^k$, and the right hand side consists of the single factor of the form ${1}/(1-q^m),$ with $m=1$. 

Recall the proof of this formula. 
We first prove the finite version, the sum of the first $n$ terms of a geometric progression:
$$1+q+q^2+\cdots+q^{n-1}=\frac{1-q^n}{1-q}.$$
Let $$S_n=1+q+q^2+\cdots + q^{n-1}.$$
Multiply both sides by $(1-q)$ to obtain:
\begin{eqnarray*}
S_n(1-q)&=&(1-q)+(q-q^2)+(q^2-q^3)+\cdots +(q^{n-1}-q^{n})\\
&=& 1-q^{n}\\
\implies S_n &=& \frac{1-q^{n}}{1-q}.
\end{eqnarray*}
To go to the infinite series, we first write the sum $S_n$ as
$$S_n=\frac{1}{1-q}-\frac{q^n}{1-q}$$ 
and then take the limit as $n\to \infty$. Note that if $q$ is a real number between $-1$ and $1$, then $q^n\to 0$ as $n\to \infty$. Thus we can say that
$$\lim_{n\to\infty} S_n = \frac{1}{1-q}.$$
For those of you who are familiar with complex numbers, we may take $q$ to be a complex number that satisfies $|q|<1$. In either case, we can say that the infinite series on the LHS of \eqref{geometric-series}  {\bf converges} to the expression on the RHS.  This completes the proof of the formula for the sum of the geometric series.

To summarize, one way to make sense of the infinite series is to first consider finite sums and then take limits. 
There is a similar approach for infinite products which we will not get into. 


There is another approach to  infinite series and products that can help us avoid questions of convergence.  In this {\em formal} approach, we consider a {\bf formal power series}, an expression of the form 
$$\sum_{k=0}^{\infty} a_k q^k := a_0+a_1q+a_2q^2+a_3 q^3+\cdots,$$
where we have used the $\Sigma$ notation for writing series in shorthand. We say that the {\bf coefficient} of $q^k$
is $a_k$. 
We regard this infinite series as an \lq infinite degree polynomial' which we can add and subtract (and even multiply) term-wise, just by extending the rules of algebra of polynomials. Two series are considered equal, if for each $k$, the coefficient of $q^k$ is the same in both the series. Here we regard $q$ as an {\bf indeterminate}, and use it as a symbol without giving it any value! 

In the formal approach, we can use the geometric series \eqref{geometric-series} to expand each term  of the form $1/(1-q^m)$ into the infinite series  
$$1+q^m+q^{2m}+\cdots$$
and then multiply them together, to write sums and products as a formal power series. Of course, we can only multiply a finite number of terms at one time. 

For example, the first few terms of the LHS of the first Rogers--Ramunujan identity \eqref{rr1-expanded} can be written as
$$1+q(1+q+q^2+q^3+q^4+\cdots)+q^4(1+q+q^2+\cdots )(1+q^2+q^4+q^6+\cdots )+\cdots.$$
We can see that the first few terms are  $$1+q+q^2+q^3+2q^4+\cdots.$$
Try this approach to find the first few terms of the RHS of the identity \eqref{rr1-expanded}  and check that they match. 
(I heartily recommend that you calculate more terms on both sides, using a  computer algebra package such as Maxima, Maple or Mathematica. Maxima is freely available on the Internet. Download and use!)

To summarize, in the formal approach to understand infinite identities, we regard both sides to be formal power series, without regard to whether the series converge.  We regard the identity to be true if for each  $k$, the coefficients of $q^k$ match on both sides.  But to be able to calculate any coefficient, the computation should involve only a finite number of algebraic operations. 

For the purpose of this article, we ignore issues of convergence altogether, and only consider these identities at a formal level. 

%
%
%
%
%
%

\section{The Simplest Continued Fraction}\label{sec:simplest}
There is one more interesting infinite object that we will briefly examine before we begin our study of Ramanujan's formulas. Consider the continued fraction
\begin{equation*}\label{golden-mean-cfrac}
1+\frac{1}{1+\cfrac{1}{1+\cfrac{1}{1+\cdots}}}.
\end{equation*}
To analyze this continued fraction, we truncate the continued fraction at various places and compute the values of the resulting (finite) continued fractions. These truncated fractions are called {\bf convergents}. 
\begin{eqnarray*}
1 =& 1 = \frac{1}{1}; \cr
1+\frac{1}{1} =& 2 = \frac{2}{1}; \cr\cr
1+\frac{1}{1+\cfrac{1}{1}}=&1+\frac{1}{2}=\frac{3}{2};\cr \cr\cr
1+\frac{1}{1+\cfrac{1}{1+\cfrac{1}{1}}}=&1+\frac{1}{3/2}=\frac{5}{3};\cr\cr\cr
1+\frac{1}{1+\cfrac{1}{1+\cfrac{1}{1+\cfrac{1}{1}}}}=&1+\frac{1}{5/3}=\frac{8}{5}.
\end{eqnarray*}
The convergents are:
$$\frac{1}{1}, \frac{2}{1}, \frac{3}{2}, \frac{5}{3}, \frac{8}{5}, \dots .$$
Note that to calculate a particular fraction, we use the value of the previous convergent, so the $n$th convergent $w_n$ satisfies the following {\bf recurrence relation}:
$$w_n= 1+\frac{1}{w_{n-1}}.$$
Further, note that the denominator of the $n$th convergent is the numerator of the previous convergent.
This suggests the substitution
$$w_n=\frac{F_n}{F_{n-1}}.$$
Substituting this value in the recurrence relation yields:
$$\frac{F_n}{F_{n-1}} = 1+ \frac{1}{F_{n-1}/F_{n-2}},$$
or
$$F_n = F_{n-1}+F_{n-2}$$
after clearing out the denominators. 

Let $F_0=1$ and $F_1=1$. Then we generate the rest of the sequence using this rule to obtain
$$1, 1, 2, 3, 5, 8, 13, 21, \dots $$
The convergents are easy to determine too, by taking ratios of successive numbers in this sequence. This sequence,  known as the {\bf Fibonacci Sequence},  starred prominently in the novel (and movie) {\em The da Vinci Code}. 

I suggest that you compute ratios $F_n/F_{n-1}$ for $n =1, 2, 3, \dots $ and convince yourself that  the ratios converge. Use a computer. Can you figure out the {\em exact} value of the limit? The limit (to which the continued fraction converges) is called the {\bf golden mean}.

If you find continued fractions attractive, then you will love this Olds \cite{olds} book. Continued fractions have appeared earlier in Resonance: Shirali \cite{shirali-2000} gives a clear introduction, and gives examples of such fractions related to the number $e$; Sury \cite{sury-2005} explains their appearance in the context of special functions. 

%
%

\section{The Rogers--Ramanujan Continued Fraction}
We are now ready to dive into Ramanujan's world. The first step is to generalize the continued fraction by adding an additional parameter $q$ to it. 

The $q$-generalization of 
$$1+1+1+\cdots+1=n$$
is the familiar sum of the Geometric Progression
$$1+q+q^2+\cdots+q^{n-1}=\frac{1-q^n}{1-q}.$$
Similarly, one would like to generalize the continued fraction
\begin{equation*}\label{golden-mean-cfrac}
1+\frac{1}{1+\cfrac{1}{1+\cfrac{1}{1+\cdots}}}
\end{equation*}
to
\begin{equation}\label{RRcfrac1}
1+\frac{q}{1+\cfrac{q^2}{1+\cfrac{q^3}{1+\cdots}}}.
\end{equation}

To gain some intuition, truncate the continued fraction at different places and evaluate the resulting (finite) fractions. 
\begin{eqnarray*}
1 &=& 1;\cr \cr
1+\frac{q}{1} &=& 1+q; \cr\cr
1+\frac{q}{1+\cfrac{q^2}{1}}&=&1+\frac{q}{1+q^2}=\frac{1+q+q^2}{1+q^2}; \cr\cr\cr
1+\frac{q}{1+\cfrac{q^2}{1+\cfrac{q^3}{1}}}&=&1+\frac{q(1+q^3)}{1+q^2+q^3}
=\frac{1+q+q^2+q^3+q^4}{1+q^2+q^3}.
\end{eqnarray*}
At each step, we would like to use the previous calculation, just as in \S \ref{sec:simplest}. 
For example, to compute
$$1+\frac{q}{1+\cfrac{q^2}{1+\cfrac{q^3}{1}}}$$
we would like to use
$$1+\frac{q}{1+\cfrac{q^2}{1}}.$$
But the powers of $q$ are not quite right.  However, if we have an additional parameter $z$
$$1+\frac{zq}{1+\cfrac{zq^2}{1}},$$
then we can adjust the powers of $q$, by replacing $z$ by $zq$. 

So consider:
\begin{equation}\label{RRcfrac2}
c(z,q)=1+\frac{zq}{1+\cfrac{zq^2}{1+\cfrac{zq^3}{1+\cdots}}}.
\end{equation}
That gives us a recurrence relation:
\begin{equation}\label{rcfrac1}
c(z,q)=1+\frac{zq}{c(zq,q)}.
\end{equation}
The continued fraction $c(z,q)$ is called the Rogers--Ramunujan continued fraction.

\section{The sum side}
The sum sides of the Rogers--Ramunujan identities appear as we solve the recurrence relation \eqref{rcfrac1}. 

To gain further insight, we again turn to computation. The first few fractions are as follows: 
\begin{eqnarray*}
c_0(z,q)&=&1 = 1;\cr \cr
c_1(z,q)&=&1+\frac{zq}{1} = 1+zq;\cr\cr
c_2(z,q)&=&1+\frac{zq}{1+\cfrac{zq^2}{1}}=1+\frac{zq}{1+zq^2}=\frac{1+zq+zq^2}{1+zq^2};\cr\cr
c_3(z,q)&=&1+\frac{zq}{1+\cfrac{zq^2}{1+\cfrac{zq^3}{1}}}=1+\frac{zq(1+zq^3)}{1+zq^2+zq^3}\cr\cr
&=&\frac{1+zq+zq^2+zq^3+z^2q^4}{1+zq^2+zq^3};\cr\cr
c_4(z,q)&=&1+\frac{zq}{1+\cfrac{zq^2}{1+\cfrac{zq^3}{1+\cfrac{zq^4}{1}}}}=
1+\frac{zq(1+zq^3+zq^4)}{1+zq^2+zq^3+zq^4+z^2q^6} \cr\cr
&=&\frac{1+zq+zq^2+zq^3+zq^4+z^2q^4+z^2q^5+z^2q^6}{1+zq^2+zq^3+zq^4+z^2q^6}.
\end{eqnarray*}
The pattern is clear. 

Let $H_n(z,q)$ be the numerator of $c_n(z,q)$. 
Its denominator is $H_{n-1}(zq,q)$.
Indeed, for $n =2, 3, 4$, we have
$$c_n(z,q)=\frac{H_n(z,q)}{H_{n-1}(zq,q)} = 1+\frac{zqH_{n-2}(zq^2,q)}{H_{n-1}(zq,q)}.$$
This suggests the substitution:
$$c(z,q)=\frac{H(z,q)}{H(zq,q)}.$$
Substituting for $c(z,q)$ in the recurrence relation \eqref{rcfrac1} yields:
\begin{equation}\label{rcfrac2}
H(z,q)=H(zq,q)+zqH(zq^2,q).
\end{equation}

So far, our calculations have been analogous to our work in \S \ref{sec:simplest}. To solve this recurrence relation we require an additional trick. We {\em assume} that the solution has the power series expansion: 
$$H(z,q)=\sum_{k=0}^{\infty} a_k z^k.$$
If the solution is really of this form, then it must satisfy the recurrence relation. We will use this idea to actually compute the coefficients $a_k$ of the solution.  

By substituting the power series expansion in \eqref{rcfrac2}, we obtain
$$\sum_{k=0}^{\infty} a_kz^k 
= \sum_{k=0}^{\infty} a_kq^kz^k
+\sum_{k=0}^{\infty} a_kq^{2k+1}z^{k+1}.
$$
On comparing coefficients of $z^k$ on both sides, we get, for $k=1, 2, 3, \dots $:  
$$a_k=a_kq^k +a_{k-1}q^{2k-1},
\text{ or } a_k=\frac{q^{2k-1}}{1-q^k}a_{k-1}.$$
On iteration, this yields: 
$$a_k=\frac{q^{k^2}}{(1-q)(1-q^2)\cdots (1-q^k)}a_0.$$
Take $a_0=1$ to find that
$$H(z,q)=\sum_{k=0}^{\infty} \frac{q^{k^2}}{(1-q)(1-q^2)\cdots (1-q^k)}z^k
=\sum_{k=0}^{\infty} \frac{q^{k^2}}{\qrfac{q}{k}}z^k,$$
where to shorten our displays, we have used the notation:
$$\qrfac{q}{k}:=(1-q)(1-q^2)\cdots (1-q^k).$$

We are interested in \eqref{RRcfrac1}, which is $c(1,q)$:
\begin{equation*}
c(1,q)
=
\frac{H(1,q)}{H(q,q)};
\end{equation*}
or, in other words,
\begin{equation}
1+\frac{q}{1+\cfrac{q^2}{1+\cfrac{q^3}{1+\cdots}}}
=
\frac{\displaystyle\sum_{k=0}^{\infty} \frac{q^{k^2}}{\qrfac{q}{k}}}
{\displaystyle\sum_{k=0}^{\infty} \frac{q^{k^2+k}}{\qrfac{q}{k}}}.
\end{equation}
Thus, the Rogers--Ramanujan continued fraction is the ratio of the sums:
$$\sum_{k=0}^{\infty} \frac{q^{k^2}}{\qrfac{q}{k}}=1+\frac{q}{(1-q)}+\frac{q^4}{(1-q)(1-q^2)}
+\frac{q^9}{(1-q)(1-q^2)(1-q^3)}+\cdots$$
and 
$$\sum_{k=0}^{\infty} \frac{q^{k^2+k}}{\qrfac{q}{k}}=1+\frac{q^2}{(1-q)}+\frac{q^6}{(1-q)(1-q^2)}
+\frac{q^{12}}{(1-q)(1-q^2)(1-q^3)}+\cdots.$$
Here they are---the sum sides of the Rogers--Ramanujan identities.  

\section{The product side}\label{sec:RRProdmotivate}
The sum sides of the Rogers--Ramunujan identities arose from the analysis of the Rogers--Ramanujan continued fraction. Next,  we show how one can guess the product sides. 
Let us first consider the LHS of the first Rogers--Ramunujan identity \eqref{rr1-expanded}:
$$H(1,q)=1+\frac{q}{(1-q)}+\frac{q^4}{(1-q)(1-q^2)}
+\frac{q^9}{(1-q)(1-q^2)(1-q^3)}+\cdots.$$


Note that, in view of \eqref{geometric-series},  the sum $H(1,q)$ is of the form
$$1+q+\text{ higher powers of $q$}.$$
There is a useful to  trick to eliminate the smallest power of $q$ (in this case $q^1$) on the RHS: Multiply both sides by $1-q$.
\begin{align*}
H(1,q)(1-q) & = 1+\frac{q^4}{1-q^2}+\frac{q^9}{(1-q^2)(1-q^3)}+\cdots\\
&=1+ q^4 + \text{ higher powers of $q$}.
\end{align*}
Now the smallest power on the RHS is $q^4$. So multiply both sides by $1-q^4$ to obtain:
\begin{align*}
H(1,q)(1-q)(1-q^4) & = 1+q^6+\frac{q^9(1+q^2)}{1-q^3}+\frac{q^{16}(1+q^2)}{(1-q^3)(1-q^4)}+\cdots\\
&=1+ q^6 + \text{ higher powers of $q$}.
\end{align*}
Next, multiply by $1-q^6$ to obtain:
\begin{align*}
H(1,q)(1-q)(1-q^4)(1-q^6) & = 1+q^{9}+q^{11}+q^{14}+
\frac{q^{16}(1+q^2)(1+q^3)}{(1-q^4)}+\cdots\\
&=1+ q^9 + \text{ higher powers of $q$}.
\end{align*}
We can continue in this fashion  and hope that
$$H(1, q)(1-q)(1-q^4)(1-q^6)(1-q^9)(1-q^{11})(1-q^{14})(1-q^{16})\cdots =1.$$
Note the powers 
$1, 4, 6, 9, 11, 14, 16$. 
The guess is that the powers are $5m+1$ and $5m+4$, for  $m=0,1,2, 3, \dots$. Thus we have the conjecture:
\begin{equation}\label{rr1}
H(1, q)=\sum_{k=0}^{\infty} \frac{q^{k^2}}{\qrfac{q}{k}}
=\prod_{m=0}^{\infty}
\frac{1}{(1-q^{5m+1})(1-q^{5m+4})},
\end{equation}
where we have used the $\Pi$ notation for products, which is analogous to  the $\Sigma$ notation for sums.

The reader can check that a similar approach on $H(q,q)$ yields
$$H(q, q)(1-q^2)(1-q^3)(1-q^7)(1-q^8)(1-q^{12})(1-q^{13})(1-q^{17})\cdots =1.$$
So we can conjecture that:
 \begin{equation}\label{rr2}
 H(q, q)=\sum_{k=0}^{\infty} \frac{q^{k^2+k}}{\qrfac{q}{k}}
=\prod_{m=0}^{\infty}
\frac{1}{(1-q^{5m+2})(1-q^{5m+3})}.
\end{equation}
Equations \eqref{rr1} and \eqref{rr2}  are the Rogers--Ramanujan Identities.

The trick used here works in other contexts too. Try it on
the Riemann zeta function
$$\zeta(s) = 1+\frac{1}{2^s}+\frac{1}{3^s} + \frac{1}{4^s}+\frac{1}{5^s}+\cdots$$
and discover for yourself a famous identity of Euler!

\section{Credits and closing remarks}
We have seen an approach to discover the Rogers--Ramunujan identities. First of all, we gave a $q$-generalization of the continued fraction for the golden mean to obtain the Rogers--Ramunujan continued fraction. Next, an analysis of this continued fraction led to the sum sides of the two Rogers--Ramunujan identities. This part of our approach is based on remarks by Askey \cite{askey1}, who suggested that Ramanujan could have discovered the Rogers--Ramanujan identities in this manner. Our method to  conjecture the product sides in \S \ref{sec:RRProdmotivate} is the same as the one used by Euler \cite{euler-thesis} to find a product representation of the Riemann zeta function.

Our computations do not amount to a proof of the Rogers--Ramunujan identities. If you want  a proof, look 
at Andrews \cite{and10lectures} or
Andrews and Baxter \cite{andbax}. A longer introduction to Ramanujan's mathematics has been given by Berndt \cite{berndt1}, a book I heartily recommend to you. See Sury \cite{sury-review} for an informative book review of this book. 

Meanwhile, I hope this article has given you, dear reader, a small taste of Ramanujan's own approach to mathematics---described by Hardy \cite[p.~xxxv]{hardy} as follows:
\begin{quote}
He worked, far more than the majority of modern mathematicians, by induction from numerical examples; 
\dots
But with his memory, his patience, and his power of calculation, he combined a power of generalisation, a feeling of form, and a capacity for rapid modification of his hypotheses, that were often really startling, and made him, in his own peculiar field, without a rival on his day.
\end{quote}

%
%

\end{document}